\documentclass[a4paper,11pt]{article}
 \usepackage[utf8]{inputenc}

\usepackage{array}

 \usepackage{amsthm,amsmath,amsxtra,amscd,amssymb,xypic,color,mathtools,xspace}
\usepackage{todonotes}
\usepackage{enumerate}
\usepackage{multirow}
\usepackage{tikz}
\usepackage{hyperref}

\usepackage{ucs}
\usepackage[utf8]{inputenc}
\usepackage[german,english,spanish,french]{babel}
\usepackage[T1]{fontenc}

\usepackage{caption}
\usepackage{subcaption}

 \usepackage{gfsartemisia}
\usepackage{amssymb,amsfonts,amsmath,amsthm,mathrsfs}
\usepackage{paralist}
\usepackage{epsfig}
\usepackage{psfrag}
\usepackage[all]{xy}
\usepackage{graphicx}
\usepackage{color}

\usetikzlibrary{calc,matrix,arrows,shapes,decorations.pathmorphing,decorations.markings,decorations.pathreplacing}
\usepackage{xifthen}
\usepackage{verbatim}
\usepackage{xspace}

\usepackage{csquotes}
\usepackage[]{appendix}

 \newcommand{\NN}{{\mathbb N}} 
\newcommand{\ZZ}{{\mathbb Z}} 

\newcommand{\RR}{{\mathbb R}} 
\newcommand{\CC}{{\mathbb C}}

 \newcommand{\PP}{{\mathbb P}}
 \newcommand{\calO}{{\mathcal O}}
 
 \newcommand{\pgcd}{{\rm pgcd}}
 
 \newcommand{\moduli}[1][g]{{\mathcal M}_{#1}}
\newcommand{\omoduli}[1][g]{{\Omega\mathcal M}_{#1}}
 \newcommand{\komoduli}[1][g]{{\Omega^{k}\mathcal M}_{#1}}
 
\theoremstyle{plain}{
\newtheorem{theorem}{Th\'eor\`eme}[section]
\newtheorem{cor}[theorem]{Corollaire}

\newtheorem{prop}[theorem]{Proposition}
\newtheorem{lem}[theorem]{Lemme}

\newtheorem{thm}{Th\'eor\`eme}

}

 \theoremstyle{remark}{
\newtheorem{rem}[theorem]{Remarque}

}

\theoremstyle{definition}{
\newtheorem{defn}[theorem]{D\'efinition}

}

\newcommand{\Weierstrass}{Weierstra\ss\xspace}
 
 \hypersetup{pdfborder=0 0 0, 
	    colorlinks=true,
	    citecolor=purple,
	    linkcolor=violet,
	    urlcolor=black,
	    pdfauthor={Quentin Gendron}
	   }

 \title{\'Equation de Pell-Abel et applications}
 \author{Quentin Gendron}

\begin{document}

 \maketitle

 \selectlanguage{english}

\begin{abstract}
In this paper, we show that there are solutions of every degree $r$ of the equation of Pell-Abel on some  real hyperelliptic curve of genus $g$ if and only if $ r > g$. This result, which is known to the experts, has consequences, which seem to be unknown to the experts. First, we deduce the existence of a primitive $k$-differential on an hyperelliptic curve of genus $g$ with a unique zero of order $k(2g-2)$ for every $(k,g)\neq(2,2)$. Moreover, we show that there exists a non Weierstrass point of order~$n$ modulo a Weierstrass point on a hyperelliptic curve of genus $g$ if and only if $n > 2g$.
\end{abstract}

\selectlanguage{french}

\begin{abstract}
Dans cet article, je montrerai que l'équation de Pell-Abel possède une solution de degré $r$ sur une courbe hyperelliptique de genre $g$ si et seulement si $r >g$. Ce résultat est connu des experts, mais il possède des conséquences intéressantes, qui ne semblent pas connues. Tout d'abord, cela implique l'existence d'une $k$-différentielle primitive avec un unique zéro d'ordre $k(2g-2)$ sur une courbe hyperelliptique de genre $g$ pour tout  $(k,g)\neq(2,2)$. De plus, cela implique qu'il existe des courbes hyperelliptiques avec un point qui n'est pas de Weierstrass et qui est d'ordre $n$ modulo un point de Weierstrass si et seulement si $n >2g$. 
\end{abstract}

\setcounter{tocdepth}{1}
\tableofcontents

\bigskip
 \section{Introduction}
 
 Soit $D$ un polynôme de degré pair $2g+2$ dont les zéros sont simples. On lui associe une courbe hyperelliptique $C_{D}$ de genre~$g$ par l'équation $y^{2} = D(x)$. L'{\em équation de Pell-Abel} sur $C_{D}$ est l'équation polynomiale
 \begin{equation}\label{eq:Pellintro}
     P^2 - DQ^2 = c   \,,
  \end{equation}
 avec $ c \in \CC^\ast$. Une solution $(P,Q)$ de cette équation est dite de {\em degré} $r$ si le degré de $P$ est égal à $r$. De plus, cette solution est {\em primitive} si le degré de $P$ est minimal parmi les solutions distinctes de $(\pm\sqrt{c},0)$. 
 
 Cette célèbre équation a été étudiée en profondeur  dans différents contextes depuis son introduction par Abel \cite{Abelpell}. Dans son article, Abel donne des conditions nécessaires et suffisantes sur le polynôme $D$ pour que l'équation de Pell-Abel ait une solution. Cette équation a été par la suite réétudiée par Tchebychev \cite{cheby} en lien avec les polynômes qui s'écartent le moins possible de l'axe des abscisses sur une union d'intervalles. Par la suite, cette étude a été approfondie par Achieser, Markoff, Zolotareff et beaucoup d'autres. On pourra par exemple consulter les textes \cite{bopoly,SerreBour,soyu} et les références de ces textes pour plus de détails. Sans prétendre à l'exhaustivité, il est intéressant de noter que cette équation est aussi apparue en lien avec la physique mathématique \cite{buzhabel},  les systèmes dynamiques \cite{mcmtor}, les billards dans les ellipsoïdes \cite{dragolivre} ou les  points d'ordres finis sur les courbes hyperelliptiques définies sur le corps des rationnels  \cite{pakocombi}.
 
 \smallskip
 \par
 \paragraph{\bf Solutions de Pell-Abel.}
 Dans cet article je me propose tout d'abord d'étudier les degrés des solutions de l'équation de Pell-Abel sur certaines courbes hyperelliptiques. Une courbe hyperelliptique $C_{D}$ est {\em totalement réelle} si les racines de $D$ sont contenues dans~$\RR$.  Je vais en particulier montrer le résultat suivant. 
 \begin{thm}\label{thm:principal}
  Pour tout $g\geq 0$, il existe une courbe hyperelliptique totalement réelle $C_{D}$ de genre $g$ telle que l'équation de Pell-Abel possède une solution primitive de degré $r$  si et seulement si $r > g$. 
 \end{thm}

 De fait, je montrerai que l'équation de Pell-Abel sur une courbe hyperelliptique quelconque n'a jamais de solutions pour $r\leq g$. Pour construire des solutions pour tout $r>g$, j'utiliserai la théorie des différentielles abéliennes sur les courbes algébriques complexes (ou surfaces de Riemann). En effet, on peut associer une différentielle canonique à chaque courbe hyperelliptique définie par un polynôme~$D$ comme ci-dessus. Je monterai que cette différentielle est étroitement liée à l'équation de Pell-Abel, ce qui me permettra de montrer le résultat ci-dessus.
 
 Ce résultat est sans aucun doute connu des spécialistes. Des versions plus ou moins explicites de celui-ci peuvent se trouver dans \cite{soyu}, \cite{bopoly} et surtout \cite{bogaeff}. Toutefois, comme l'écrit Serre \cite[paragraphe 2.5]{SerreBour} à propos de résultats similaires: "Il n'est pas facile de dire à qui ces théorèmes sont dus. L'une des difficultés est qu'ils sont rarement énoncés explicitement, et, du coup, ils ne sont pas démontrés en détails".  Aussi je pense qu'il est intéressant d'énoncer  clairement ce résultat important et de le prouver de manière relativement élémentaire.  
\smallskip
\par
 Je vais maintenant donner deux applications intéressantes de ce résultat que je n'ai pu trouver dans la littérature. Les sujets que je vais aborder reflètent mes goûts et il est probable que des applications similaires à d'autres problèmes sont possibles.
 
  \smallskip
 \par
 \paragraph{\bf Les $k$-différentielles.}
 Je voudrais tout d'abord utiliser le théorème~\ref{thm:principal} dans le domaine des $k$-différentielles sur les surfaces de Riemann.  Soit $k\geq1$, une {\em $k$-différentielle} est une section non nulle de la puissance tensoriel $k$ième du fibré canonique d'une surface de Riemann. Une $k$-différentielle est {\em primitive} si elle n'est pas la puissance d'une $k'$-différentielle avec $0 < k' < k$.
 
 Un problème délicat est de savoir s'il existe des $k$-différentielles satisfaisant certaines propriétés. Dans cet direction, je montrerai le résultat suivant.
 \begin{thm}\label{thm:proj}
Soit  $g\geq2$ et $k\geq1$, il existe une  $k$-différentielle primitive avec un unique zéro d'ordre $k(2g-2)$ sur une courbe hyperelliptique totalement réelle de genre $g$ si et seulement si $(g,k)\neq(2,2)$.
 \end{thm}
 Notons que le théorème 1.4 de \cite{geta} donne l'existence de $k$-différentielles primitives avec un unique zéro sur des surfaces de Riemann de genre $g$, sauf dans le cas où $g=k=2$. L'originalité du présent résultat est donc de donner l'existence de tels objets sur une courbe hyperelliptique. Le théorème~\ref{thm:proj} est une application directe du théorème~\ref{thm:principal} et de résultats élémentaires de géométrie algébrique.
 \par
 De plus, ce théorème est un cas particulier d'un résultat plus général. On peut en effet se demander s'il existe des $k$-différentielles holomorphes primitives sur des courbes hyperelliptiques dont le diviseur est supporté sur des points conjugués par l’involution hyperelliptique. J’étudie ce problème dans la section~\ref{sec:kdiffconjug}. En particulier, la proposition~\ref{prop:kdiffconjug} donne des conditions nécessaires et suffisantes pour l'existence de telles $k$-différentielles.
 \smallskip
 \par
 \paragraph{\bf Les points de torsion.}
 La seconde application que je souhaite donner est au sujet  des points de torsion sur les courbes hyperelliptiques.   
 Étant donnés une courbe  hyperelliptique~$C$ et point de \Weierstrass~$W$ de $C$, on considère l'application d'Abel-Jacobi $\varphi_{W}\colon C \to J_{C}$ qui envoie le point  $W$  sur l'origine de la jacobienne $J_{C}$ de~$C$. Un {\em point de $n$-torsion primitif} modulo~$W$ sur $C$  est un point $P$ tel que $\varphi_{W}(P)$ est un point de $n$-torsion de la jacobienne qui n'est pas de $n'$-torsion pour tout $0< n' < n$.

 \begin{thm}\label{thm:torsion}
Soit $g\geq 2$ et $n > 2$, il existe une courbe hyperelliptique totalement réelle de genre $g$, un point de \Weierstrass $W$ sur $C$  et un point $P\in C$, tels que $P$ est un point de $n$-torsion si et seulement si $n > 2g$.
 \end{thm} 
 De  fait, je montrerai qu'il n'existe aucune courbe hyperelliptique de genre $g$ avec un point de $n$-torsion primitif pour $2<n<2g$.
 De nombreux exemples de points de $n$-torsion sont déjà connus (voir par exemple \cite{beza,lepre} entre beaucoup d'autres) et le théorème~\ref{thm:torsion} est prouvé en genre $2$ dans le paragraphe~4 de~\cite{bogrtheta}. Toutefois, ce résultat ne semble pas connu pour $g\geq3$. La preuve de celui-ci se fait par récurrence sur le genre des courbes hyperelliptiques, en associant les techniques de dégénération des différentielles abéliennes de \cite{BCGGM1} et la relation entre la différentielle canonique et l'équation de Pell-Abel détaillée dans la section~\ref{sec:pellabel}.

 \smallskip
 \par
 \paragraph{\bf Structure.}
Cet article est organisé de la façon suivante. Dans la section~\ref{sec:pellabel}, je donne les rappels sur l'équation de Pell-Abel et la preuve du théorème~\ref{thm:principal} sur le degré des solutions. La section~\ref{sec:kdiff} est dédiée à l'application aux $k$-différentielles. Le théorème~\ref{thm:proj} et une généralisation de celui-ci y sont prouvés. Dans cette section, je donne aussi quelques rappels sur les $k$-différentielles, qui peuvent être utiles à la section précédente, et sur leurs dégénérations. Dans la section~\ref{sec:tor}, après quelques rappels sur les points de $n$-torsion, je prouverai le théorème~\ref{thm:torsion}.

 \smallskip
 \par
 \paragraph{\bf Remerciements.}
 Cet article vient du mon envie de comprendre la partie liée aux différentielles abéliennes du séminaire Bourbaki de Jean-Pierre Serre~\cite{SerreBour}. Je remercie donc Nicolas Bourbaki d'organiser son séminaire et de mettre les vidéos en accès libre sur internet. Je remercie   Mikhail Sodin de m'avoir procuré une copie de la traduction anglaise de \cite{soyu} et John Boxall  pour avoir porté à ma connaissance certaines références. Enfin je remercie Jean-Pierre Serre et Andrei Bogatyrëv pour leurs commentaires bienveillants sur un brouillon de ce texte.

 \section{Équation de Pell-Abel et différentielle canonique}
 \label{sec:pellabel}
 
 Le but de cette section est de prouver le théorème~\ref{thm:principal} sur les solutions de l'équation de Pell-Abel. Pour cela, je commencerai  par une discussion de cette équation et de la différentielle canonique sur les courbes hyperelliptiques totalement réelles. Elle s'inspire du  paragraphe~2 de \cite{SerreBour} dont je  suivrai les notations.  D'autres bonnes références sur l'équation de Pell-Abel sont  \cite{bopoly} et \cite{soyu}, où de nombreuses autres sources sont données et \cite{bogaeff}, dont la section~2 est proche de cette section.

 \subsection{L'équation de Pell-Abel sur les courbes hyperelliptiques}\label{sec:pelldebut}
 On fixe un entier $g\geq 0$ et on se donne $2g+2$  nombres réels  
 $$a_0 < b_0 < a_1 < \cdots < a_g < b_g\,.$$ 
 Pour chaque $j = 0, \dots, g$ on note $E_j = [a_j,b_j]$ et $E$ la réunion des $E_j$. On définit le polynôme
  $$D(x) = \prod_{i=0}^{g}(x-a_{i})(x-b_{i})$$
  et  
  la courbe hyperelliptique affine~$C^{\rm aff}_D$ de genre  $g$ donnée par l'équation $y^2 = D(x)$. Sa complétée $C_D$ s'obtient en lui ajoutant deux points à l'infini, que l'on note $\infty_+$ et~$\infty_-$. Ces points sont caractérisés par le fait que $y/x^{g+1}$ prend respectivement la valeur  $1$  au premier et $-1$ au second. Une courbe hyperelliptique ramifiée uniquement au dessus de point réels, comme $C_{D}$, est appelée {\em totalement réelle}.   Notons que  $D(x) $ est inférieur ou égal à $ 0$ si $x \in E$ et  strictement supérieur à $ 0$ sinon. 
 
  L'algèbre affine de $C^{\rm aff}_D$ est un $k[x]$-module libre de base $\{1,y\}$. Ses éléments inversibles
  sont de la forme:
  \begin{equation}\label{eq:Pell}
     f = P + yQ,  \text{ avec }P,Q \in \CC[x]  \text{ et }  P^2 - DQ^2 = c  \text{ où } c \in \CC^\ast.
  \end{equation}
 On supposera que les polynômes $P$ et $Q$ sont unitaires et que $Q$ est  non nul. Ce second point revient à ne pas considérer les solutions triviales $(\pm\sqrt{c},0)$.  L'équation  $P^2 - DQ^2 = c$ est appelée {\em équation de Pell-Abel} car elle apparaît pour la première fois dans \cite{Abelpell} généralisant l'équation diophantienne de Pell. Le {\em degré} d'une solution $f$ est le degré de~$P$. On dira que~$f$ est une solution {\em primitive} si le degré de $P$ est minimal (parmi les solutions distinctes de $f = \pm\sqrt{c}$). 
 
  Une différence importante entre les équations de Pell celle  de Pell-Abel est que cette dernière n'a pas toujours de solutions non triviales.
 De façon plus précise, soit $r$ un entier supérieur ou égal à~$1$. Il y a équivalence entre les deux conditions suivantes.
 \begin{itemize}
  \item[(i)] L'équation~\eqref{eq:Pell} possède une solution de degré $r$.
  \item[(ii)] Le diviseur $r(\infty_+ - \infty_-)$ de la courbe $C_D$ est linéairement équivalent à~$0$.
 \end{itemize}
 On peut reformuler le point (ii) de la manière suivante. Soit $J_{D}$ la jacobienne de $C_D$. Le diviseur $\infty_+ - \infty_-$ définit un
  point $P_\infty$ de $J_{D}$. Le point (ii) signifie que  $rP_\infty = 0$, autrement dit que $P_\infty$ est un point d'ordre fini divisant~$r$. De plus, cet ordre est égal à $r$ si la solution est primitive de degré~$r$.
  
En genre $g=0$, le fait que $J_{C}=0$ implique donc que l'équation de Pell-Abel a des solutions pour tout $r\geq1$. Les solutions sont données par les polynômes de Tchebychev usuels  $P,Q$ respectivement  de première et seconde espèce (modulo normalisation). 
Dans le cas $g=1$, la jacobienne est égale à la courbe et l'existence de points de $r$-torsion pour $r\geq2$ est claire. 
  Si $g\geq 2$, il est facile de construire des exemples de courbes hyperelliptiques dont l'équation de Pell-Abel n'a pas de solution. Par exemple, dans le cas $g=2$, s'il existe une solution primitive de degré  $r$ de l'équation de Pell-Abel sur~$C_{D}$, alors il existe une $r$-différentielle primitive avec un unique zéro d'ordre $2r$ en $\infty_{+}$ (voir la section~\ref{sec:kdiff} pour plus de détails). On sait par le théorème~1.1 de \cite{BCGGM3} que la dimension de l'espace des modules des $r$-différentielles primitives avec un unique zéro sur les surfaces de Riemann de genre $2$ est égal à~$3$.  Donc sa projection à $\moduli[2]$ est de dimension $2$ et les courbes hyperelliptiques qui possèdent des solutions primitives de degré $r$  à l'équation de Pell-Abel forment des diviseurs (éventuellement non irréductibles) dans~$\moduli[2]$.

  \subsection{La différentielle canonique}
 Je vais maintenant introduire la {\em différentielle canonique} sur une courbe hyperelliptique totalement réelle $C_{D}$, tout d'abord arbitraire puis telle que l'équation de Pell-Abel possède une solution. Je souhaite souligner ici que cette différentielle n'est pas canoniquement associée à une courbe hyperelliptique abstraite $C$, mais qu'elle dépend du choix du polynôme $D$ définissant la surface~$C$. Quelques rappels sur les différentielles abéliennes sont donnés dans la section~\ref{sec:kdiff} et peuvent être utiles au lecteur peu familier avec cette théorie.
 \smallskip
 \par
  \paragraph{Différentielle canonique sur une courbe hyperelliptique totalement réelle.}
  On se donne une courbe hyperelliptique totalement réelle $C_{D}$.
  Étant donné  un polynôme unitaire $A \in \RR[x]$  de degré $g$, on lui associe la forme différentielle $\eta_A = \frac{A(x)}{y} dx$   sur $C_D$. Cette forme est holomorphe ailleurs qu'en l'infini et possède  un pôle simple en $\infty_+$ et en $\infty_-$, dont les résidus sont respectivement $-1$ et $+1$.  Notons que $\eta_A $ est une  forme de troisième espèce et que changer $A$ revient à ajouter à $\eta_{A}$ une forme de première espèce, i.e. une différentielle holomorphe. Puisque l'espace vectoriel des différentielles holomorphes a pour base les formes $\frac{x^j}{y}dx$ avec $0\leq j < g$, on en déduit  qu'il existe un unique polynôme $A$ tel que 
  \begin{equation}\label{eq:2.3.2}
   \int_{b_{j-1}}^{a_j} \frac{A(x)}{\sqrt{D(x)}}\,dx = 0\text{  pour } j=1,\ldots,g\,.
  \end{equation}
  Ces intégrales sont les {\em périodes réelles} de $\eta_A$.
Dans cette formule, on note $\sqrt{\cdot}$ la racine positive.  Dans la suite, si $t$ est négatif, on définit $\sqrt{t}$ comme $i\sqrt{-t}$. On peut donc introduire l'objet central de cette section.
\begin{defn}\label{def:cano}
 Le polynôme $A$ satisfaisant l'équation~\eqref{eq:2.3.2} est le {\em polynôme canonique} de~$C_{D}$ et est noté~$R$. La forme différentielle $\eta_R $ est la
  {\em différentielle canonique} de~$C_{D}$, que l'on notera simplement~$\eta$. 
\end{defn}

  La formule~\eqref{eq:2.3.2} dit que
  dans chaque intervalle intermédiaire $T_j =  [b_{j-1},a_j]$ l'intégrale de $\eta = \frac{R(x)}{\sqrt{D(x)}}$ est nulle. Cela entraîne que $R(x)$ change de signe dans $T_j$, donc s'annule en au moins un point intérieur de $T_j$. Comme le nombre des $j$  est égal au degré de $R$, on obtient le résultat suivant.
  \begin{lem}\label{lem:2.3.4}
   Le polynôme $R(x)$ a une racine et une seule à l'intérieur de chaque $T_j$ et n'a aucune autre racine (réelle ou complexe).
  \end{lem}
 Comme les zéros de~$\eta$ sont les préimages des racines de $R$, le lemme~\ref{lem:2.3.4} implique que  la différentielle canonique~$\eta$ possède exactement $2g$ zéros simples.    On note $z_{j}$ le zéro de~$\eta$ dont la projection par l'application hyperelliptique appartient à  l'intervalle $T_{j}$ et tel que la partie réelle de l'intégrale de $\eta$ de $b_{j-1}$ à $z_{j}$ est positive. Le zéro qui est conjugué à $z_{j}$ par l'involution hyperelliptique est noté $\bar{z}_{j}$.
  
  On s'intéresse maintenant aux $g+1$ intégrales
  \begin{equation}\label{eq:2.3.6}
   \eta_j = \int_{a_j}^{b_j} \frac{R(x)}{\sqrt{D(x)}} \,dx = \int_{a_j}^{b_j} \eta, \text{ pour } j=0,1,\dots,g \,.
  \end{equation}
 Les nombres $2\eta_j$ sont appelés les {\em  périodes imaginaires} de la différentielle canonique~$\eta$.
  Notons que  comme $R$ est un polynôme réel et $D(x) \leq 0$ pour $x\in E$,  les parties réelles  des~$\eta_j$ sont nulles.  
  
  Je vais maintenant énoncer une propriété importante satisfaite par les périodes imaginaires.  Cette propriété est une conséquence du théorème des résidus et du fait que les résidus aux pôles de $\eta$ sont~$\pm1$. Les détails sont donnés à la fin de la section~2.3 de~\cite{SerreBour}.
  \begin{lem}\label{lem:2.3.7}
   Il existe des signes $\varepsilon_j \in \{-1,1\}$
  tels que l'on ait:
  \begin{equation*}
   \sum_{0\leq j \leq g}  \varepsilon_j \eta_j = i\pi\,.
  \end{equation*}
  \end{lem}
  
  J'introduis maintenant les dernières quantités qui vont permettre de paramétrer les différentielles canoniques.
 Pour $j=1,2,\dots,g $ on définit les valeurs
\begin{equation}
    \lambda_j = \int_{b_j-1}^{z_j} \eta = \int_{a_{j}}^{z_j} \eta\,.
\end{equation}
On appelle les quantités $2\lambda_{j}$ les {\em périodes plates} de $\eta$. En effet, elles correspondent à l'intégrale de $\eta$ entre les zéros $\bar{z}_{j}$ et $z_{j}$.
  \smallskip
  \par
  \paragraph{Différentielle canonique si l'équation de Pell-Abel a une solution.}
  Jusqu'à présent, la discussion ne dépendait pas du fait que l'équation de Pell-Abel possède une solution sur $C_{D}$ ou non. Dans la suite de ce paragraphe, je travaillerai sous l’hypothèse que cette équation possède une solution $f = P + yQ$ de degré $r$ sur $C_{D}$. 
  \smallskip
  \par
  Le premier résultat relie la solution de Pell-Abel à la différentielle canonique.
  \begin{lem}\label{lem:logderiv}
   Soit $C_{D}$ une courbe hyperelliptique totalement réelle telle que l'équation de Pell-Abel possède une solution $f = P + yQ$ de degré $r$. La différentielle canonique~$\eta$ de $C_{D}$ est égale à
   \begin{equation}\label{eq:logderiv}
    \eta = r \frac{df}{f}
   \end{equation}
  \end{lem}
Il est clair que $df/f$ possède deux pôles simples en $\infty_+$ et $\infty_-$ de résidus respectifs~$-r$ et $r$. Pour prouver que $df/f=r\eta$, il suffit donc de prouver que les périodes réelles de la forme $df/f$ sont nulles. Pour $j = 1,\dots,g$, on considère $f$ comme une fonction~$f(x)$ sur l'intervalle $T_j= [b_{j-1},a_{j}]$ en choisissant une détermination de $y$. Cela permet de considérer
  $f$ comme une fonction~$f(x)$ sur $T_j$. Comme cette fonction est réelle, et ne s'annule pas, on en déduit que $f(b_{j-1}) = f(a_j) = \varepsilon \sqrt{c}$, avec $\varepsilon = \pm 1$. L'intégrale de $df/f$ est donc la fonction $\log(\varepsilon f)$ qui prend les mêmes valeurs aux extrémités de $T_j$. Cela entraîne que l'intégrale sur  $T_j$ de sa dérivée est nulle et donc le lemme~\ref{lem:logderiv}.
  \smallskip 
  \par
  On peut maintenant préciser le lemme~\ref{lem:2.3.7} dans le cas où l'équation de Pell-Abel est résoluble sur  $C_{D}$. Rappelons que $\eta_{i}$ est la $i$ème période imaginaire de~$\eta$.
  \begin{prop}\label{lem:2.3.7bis}
   Soit $C_{D}$ une courbe hyperelliptique totalement réelle telle que l'équation de Pell-Abel possède une solution de degré $r$ et $\eta$ sa différentielle canonique.  Les $r_j = r |\eta_j| /\pi$ sont des nombres entiers qui satisfont l'équation
  \begin{equation}\label{eq:2.3.8bis}
   \sum_{0\leq j \leq g}  |\eta_j| = r\,.
  \end{equation}
  \end{prop}
 Fixons l'un des intervalles $E_j= [a_{j},b_{j}]$ avec $j = 0,\dots,g$. On considère $f$ comme une fonction $f(x) = P(x) + iy_1(x)Q(x)$  sur~$E_j$, où $y_{1}$ est la fonction $y_1(x) = \sqrt{-D(x)}$.
  Cette fonction est à valeurs complexes et de module $\sqrt{c}$. On peut donc l'écrire sous la forme $f(x) = \sqrt{c}e^{i\vartheta(x)}$, avec $\vartheta\colon E_j \to \RR/2\pi \ZZ$. On peut relever $\vartheta$ en une fonction continue $\theta\colon E_j \to \RR$ en lui imposant sa valeur en~$a_j$. Comme $f(a_j)=\pm \sqrt{c}$, cette valeur est de la forme
 $c_0\pi$, avec $c_0 \in \ZZ$. La valeur de~$f$ en~$b_j$ est alors $c_1\pi$ avec $c_1\in \ZZ$. On a donc:
 \begin{equation*}
  P(x) = \sqrt{c}\cos(\theta(x)) \text{ et } y_1(x)Q(x) = \sqrt{c}\sin(\theta(x)).
 \end{equation*}
 De plus, par le lemme~\ref{lem:logderiv} on a:
 \[\int_{a_j}^{b_j} \frac{\eta}{r} = \int_{a_j}^{b_j} \frac{df}{f} = i  \int_{a_j}^{b_j} d\theta = i(\theta(b_j)- \theta(a_j)) =(c_1 - c_0)i\pi .\]
 Le nombre $r_j = |c_1-c_0|$ est donc un entier strictement supérieur à $0$.
 Par le lemme~\ref{lem:2.3.7}, il existe des signes~$\pm$ tels que:
 \begin{equation}\label{eq:2.5.3}
   \sum_j  \pm r_j = r .
 \end{equation}
  Il reste donc à montrer que tous les signes sont positifs. Comme $\theta$ est strictement monotone sur $E_{j}$ entre  $c_{0}\pi$ et $c_{1}\pi$, son cosinus s'annule $|c_{1}-c_{0}|$ fois. Cela implique que  le nombre de racines de $P$ dans $E_j$ est $r_j$ et comme $P$ a au plus~$r$ racines, on a donc:
 \begin{equation}\label{eq:2.5.4}
  \sum_j r_j \leq r.
\end{equation}
 L'équation~\eqref{eq:2.3.8bis} est alors obtenue en comparant les équations~\eqref{eq:2.5.3} et~\eqref{eq:2.5.4}. 
  \smallskip
  \par
  Je vais maintenant résumer les différents résultats obtenus précédemment. Rappelons que pour tout $j =0,\dots,g$, on pose $r_j = r |\eta_j| /\pi$, où $\eta_j$ est la $j$-ième demi-période imaginaire de la différentielle canonique~$\eta$ définie par l'équation~\eqref{eq:2.3.6}.
\begin{prop}\label{prop:propdiffcancara}
 Soit $\eta$ la différentielle canonique d'une courbe hyperelliptique de genre~$g$ totalement réelle~$C_{D}$ telle que l'équation de Pell-Abel possède une solution de degré~$r$.  On a alors les propriétés suivantes.
 \begin{enumerate}[(i)]
    \item La différentielle $\eta$ possède $2g$ zéros simples et $2$ pôles simples dont les résidus sont~$\pm 1$.
    \item Les $g$ périodes réelles de $\eta$ sont nulles.
    \item Les $r_j$ sont des entiers strictement positifs de somme $r$.
   \end{enumerate}
\end{prop}
   Une différentielle canonique satisfaisant à la condition (iii) de la proposition~\ref{prop:propdiffcancara} est dite de {\em degré}~$r$. De plus, si les $r_{j}$ sont premiers entre eux on dit que la différentielle est {\em primitive}. On vérifie facilement que si  $f$ est primitive, alors la différentielle est primitive.

  Cette proposition possède déjà des conséquences intéressantes sur les solutions de l'équation de Pell-Abel. En particulier, je souhaite faire la remarque suivante, qui se trouve dans le paragraphe~5.1 de \cite{soyu}.
  \begin{rem}\label{rem:nonex}
  Le point (iii) de la proposition~\ref{prop:propdiffcancara} implique est qu'il n'existe pas de solution de degré $r\leq g$ de l'équation de Pell-Abel sur une courbe hyperelliptique totalement réelle de genre~$g$.  Ce résultat est valable pour toutes les courbes hyperelliptiques complexes comme je le montrerai dans la section~\ref{sec:thmtorsion}.
  \end{rem}
  \par

\smallskip
\par
\paragraph{La structure plate d'une différentielle canonique.}
Pour terminer la discussion des différentielles canoniques, je vais décrire la représentation plate d'une différentielle qui satisfaisait aux conditions de la proposition~\ref{prop:propdiffcancara}. Le lecteur non familier avec cette théorie pourra consulter la section~\ref{sec:rappelsdiff} et les nombreuses références sur le sujet, comme par exemple \cite{mata,zorich}. Notons que par le lemme~\ref{lm:existancecano} ci-dessous, cette différentielle est une différentielle canonique (indépendamment du choix des périodes plates $\lambda_{i})$. 

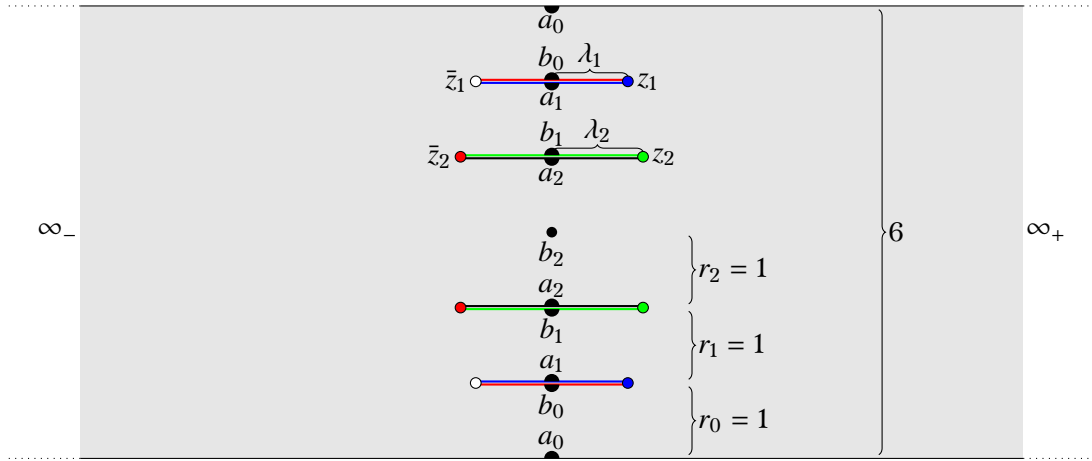
\begin{figure}[ht]
 \centering
\begin{tikzpicture}[scale=1,decoration={
    markings,
    mark=at position 0.5 with {\arrow[very thick]{>}}}]
   \fill[fill=black!10] (-6.2,-3) coordinate (r1) -- (6.2,-3) coordinate (r2) coordinate[pos=.5] (a0) -- (6.2,3) coordinate (r3) -- (-6.2,3) coordinate (r4)coordinate[pos=.5] (a0b) -- cycle;
   \draw (r1) -- (r2);
   \draw (r3) -- (r4);
   \draw[dotted] (r1) -- ++(-1,0);
    \draw[dotted] (r2) -- ++(1,0);
\draw[dotted] (r3) -- ++(1,0);
\draw[dotted] (r4) -- ++(-1,0);

\draw[red,thick] (-1,-2.02) -- (1,-2.02) coordinate[pos=.5] (b0);
\draw[blue,thick] (-1,-1.98) -- (1,-1.98) coordinate[pos=.5] (a1);
\draw[green,thick] (-1.2,-1.02) -- (1.2,-1.02) coordinate[pos=.5] (b1);
\draw[black,thick] (-1.2,-.98) -- (1.2,-.98) coordinate[pos=.5] (a2);
\draw[red,thick] (-1,2.02) -- (1,2.02) coordinate[pos=.5] (b0b);
\draw[blue,thick] (-1,1.98) -- (1,1.98) coordinate[pos=.5] (a1b);
\draw[green,thick] (-1.2,1.02) -- (1.2,1.02) coordinate[pos=.5] (b1b);
\draw[black,thick] (-1.2,.98) -- (1.2,.98) coordinate[pos=.5] (a2b);

\coordinate (b2) at (0,0);

\coordinate (c1) at (-1,-2); 
\coordinate (c2) at (-1.2,-1);
\coordinate (c3) at (-1.2,1); 
\coordinate (c4) at (-1,2);
\coordinate (c1b) at (1,-2); 
\coordinate (c2b) at (1.2,-1);
\coordinate (c3b) at (1.2,1); 
\coordinate (c4b) at (1,2);

\filldraw[fill=white] (c1) circle (2pt);
\filldraw[fill=white] (c4) circle (2pt);
\filldraw[fill=red] (c2) circle (2pt);
\filldraw[fill=red] (c3) circle (2pt);

\filldraw[fill=blue] (c1b) circle (2pt);
\filldraw[fill=blue] (c4b) circle (2pt);
\filldraw[fill=green] (c2b) circle (2pt);
\filldraw[fill=green] (c3b) circle (2pt);

\node[right] at (c3b) {$z_{2}$};
\node[right] at (c4b) {$z_{1}$};
\node[left] at (c3) {$\bar{z}_{2}$};
\node[left] at (c4) {$\bar{z}_{1}$};

\foreach \i in {0,1}{
\fill (a\i)++(.1,0) arc [start angle=0, end angle=180, radius=1mm];
\node[above] at (a\i) {$a_{\i}$}; 
\fill (a\i b)++(-.1,0)  arc [start angle=-180, end angle=0, radius=1mm];
\node[below] at (a\i b) {$a_{\i}$}; 
\fill (b\i)++(-.1,0)  arc [start angle=-180, end angle=0, radius=1mm];
\node[below] at (b\i) {$b_{\i}$}; 
\fill (b\i b)++(.1,0)  arc  [start angle=0, end angle=180, radius=1mm];
\node[above] at (b\i b) {$b_{\i}$}; 
}
\fill (a2)++(.1,0) arc [start angle=0, end angle=180, radius=1mm];
\node[above] at (a2) {$a_{2}$}; 
\fill (a2b)++(-.1,0)  arc [start angle=-180, end angle=0, radius=1mm];
\node[below] at (a2b) {$a_{2}$}; 
\fill (b2) circle (2pt);\node[below] at (b2) {$b_{2}$};

\draw [decorate,decoration={brace}]	(1.8,-2.05) -- (1.8,-2.95) node [midway, right] {$r_{0}=1$};
\draw [decorate,decoration={brace}]	(1.8,-1.05) -- (1.8,-1.95) node [midway, right] {$r_{1}=1$};	
\draw [decorate,decoration={brace}]	(1.8,-.05) -- (1.8,-.95) node [midway, right] {$r_{2}=1$};	

\draw [decorate,decoration={brace}]	(4.3,2.95) -- (4.3,-2.95) node [midway, right] {$6$};	

\draw [decorate,decoration={brace}]	(0,2.07) -- (1,2.07) node [midway, above] {$\lambda_{1}$};
\draw [decorate,decoration={brace}]	(0,1.07) -- (1.2,1.07) node [midway, above] {$\lambda_{2}$};

\node at (-6.5,0) {$\infty_{-}$};
\node at (6.5,0) {$\infty_{+}$};

\end{tikzpicture}
 \caption{Une différentielle canonique (multipliée par $\tfrac{r}{\pi}$) en genre $2$ et de degré $r=3$.} \label{fig:diffcanong2}
\end{figure}

Dans la  figure~\ref{fig:diffcanong2}, je représente la surface de translation associée à une différentielle primitive de degré~$r=3$ (après multiplication par $\pi/r$) sur une courbe hyperelliptique totalement réelle de genre~$2$. Cette surface de translation est obtenue de la façon suivante. On prend un cylindre horizontal infini de hauteur~$2r$ et on considère un point $a_{0}$ sur celui-ci. On fait $2g$ coupures horizontales de la façon suivante. Pour tout $i=1,\dots,g$ la coupure~$c_{i}^{+}$, resp. $c_{i}^{-}$, est de longueur $2\lambda_{i}$, centrée sur le point à une distance $ \sum _{j=1}^{i }r_{j}$ au dessus, resp. en dessous, de~$a_{0}$. On obtient alors deux composantes de bord pour chaque coupure, que j’appelle {\em lèvres}  supérieures et inférieures respectivement. On identifie alors par translation (verticale) la lèvre inférieure de $c_{i}^{+}$ avec la lèvre supérieure de~$c_{i}^{-}$. On identifie de même la lèvre supérieure de $c_{i}^{-}$ avec la lèvre inférieure de $c_{i}^{+}$.  Ces identifications sont représentées  sur la figure~\ref{fig:diffcanong2} par des segments de même couleur. Grâce à la correspondance exposée à la section~\ref{sec:rappelsdiff}, la surface de translation ainsi définie correspond à une différentielle abélienne.
\par
Les  points $a_{i}$ et $b_{j}$ sur cette surface de translation sont donnés de la façon suivante. Le point $b_{i}$ pour $i<g$ correspond au milieu du segment obtenu en identifiant la lèvre inférieure de $c_{i+1}^{+}$  avec la  lèvre supérieure de $c_{i+1}^{-}$. Le point $b_{g}$ est le point qui se trouve au dessus de $a_{0}$ à distance~$r$. Le point $a_{i}$ pour $i>0$ correspond au milieu du segment obtenu en identifiant la  lèvre supérieure de $c_{i}^{+}$ avec la lèvre inférieure de $c_{i}^{-}$. Enfin, les zéros $z_{i}$ et $\bar{z}_{i}$ de la différentielle correspondent respectivement aux sommets de droite et de gauche des segments $c_{i}^{\pm}$.
\par
 Je justifie maintenant que la différentielle représentée sur la figure~\ref{fig:diffcanong2} possède bien les propriétés de la proposition~\ref{prop:propdiffcancara} avec $r=3$.  Le fait que la distance entre les points~$a_{i}$ et $b_{i+1}$ est nulle implique directement le point (ii) de la proposition. Le point (iii) est clair par construction. Pour le point (i), on constate que les singularités coniques de la surface de translation sont d'angle $4\pi$. Par le dictionnaire entre surfaces de translations et différentielles, ces singularités correspondent à des zéros d'ordre~$1$ de la différentielles. Les deux fins du cylindres correspondent à deux pôles simples. Et, par le théorème des résidus, la hauteur du cylindre est égal au résidu du pôle correspondant, modulo $\pm 2i\pi$. Cela conclut la vérification du point~(i).
\par
On peut considérer la préimage de l'axe réel sur les courbes hyperelliptiques totalement réelles. Elle possède deux composantes qui s'intersectent au dessus des points de ramifications $a_{i}$ et $b_{j}$.  Il s'agit des deux demi-droites horizontale qui viennent de $\pm\infty$ et arrivent en $a_{0}$. L'image du segment  $E_{0} = [a_{0},b_{0}]$ est composé des deux segments verticaux dénotés par $\left[a_{0},b_{0}\right]$ sur la figure. Par la suite, on obtient sur la figure les segments notés $\left[b_{0},z_{1}\right]$ et $\left[b_{0},\bar{z}_{1}\right]$. Puis les deux segments  $\left[z_{1},a_{1}\right]$ et $\left[\bar{z}_{1},a_{1}\right]$. L'union de ces quatre segment est l'image du segment $T_{1}$. Puis on a les deux segments verticaux $\left[a_{1},b_{1}\right]$ et on continue ainsi jusqu'à arriver au point $b_{g}$.  Du point $b_{g}$ l'image est constituée des deux demi-droites horizontales qui pointent vers $\pm\infty$.

 \subsection{La preuve du théorème~\ref{thm:principal}}

Je vais commencer la preuve du théorème~\ref{thm:principal} par un lemme fondamental sur les différentielles canoniques. Ce résultat est similaire au théorème~2.1 de~\cite{bopoly} et au résultat énoncé dans la section~5.5 de \cite{soyu}.
\begin{lem}\label{lm:existancecano}
 Soit $\omega$ une différentielle satisfaisant aux conditions de la proposition~\ref{prop:propdiffcancara}, alors $\omega$ est une différentielle canonique de degré $r$ sur une courbe hyperelliptique totalement réelle. 
\end{lem}

\begin{figure}[ht]
 \centering
\begin{tikzpicture}[scale=1.5,decoration={
    markings,
    mark=at position 0.5 with {\arrow[very thick]{>}}}]
   \fill[fill=black!10] (0,-3) coordinate (a0) -- (6,-3) coordinate (r2)-- (6,0) coordinate (r3) -- (0,0) coordinate (a0b) -- cycle;
   \draw (r2) -- (a0) -- (a0b) -- (r3);
    \draw[dotted] (r2) -- ++(.5,0);
\draw[dotted] (r3) -- ++(.5,0);

\draw[thick] (0,-2) coordinate(b0) -- (1,-2);
\draw[thick]  (0,-1) coordinate (b1) -- (1.2,-1);

\draw [decorate,decoration={brace}]	(1.3,-2.05) -- (1.3,-2.95) node [midway, right] {$q_{1}\pi = \pi$};
\draw [decorate,decoration={brace}]	(2.7,-1.05) -- (2.7,-2.95) node [midway, right] {$q_{2}\pi = 2\pi$};	
\draw [decorate,decoration={brace}]	(4.3,-.05) -- (4.3,-2.95) node [midway, right] {$M\pi = 3\pi$};

\draw [decorate,decoration={brace}]	(0,-1.93) -- (1,-1.93) node [midway, above] {$h_{1}$};
\draw [decorate,decoration={brace}]	(0,-.93) -- (1.2,-.93) node [midway, above] {$h_{2}$};
\end{tikzpicture}
 \caption{Le peigne $\Pi_{M}(1,2;h_{1},h_{2})$ avec $M=3$.} \label{fig:peigne}
\end{figure}
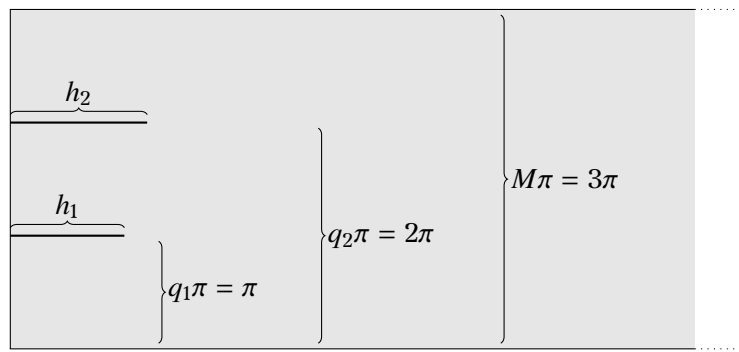

La preuve de ce résultat s'appuie sur le très bel argument d'Achieser exposé dans la section~5.4 de \cite{soyu}. Elle repose sur la formule de Schwarz-Christoffel, qui est par exemple exposée dans \cite[IV.2.3]{dSG}, et qui s'énonce de la façon suivante.  Étant donné un domaine polygonal simplement connexe~$P$  dont les  sommets sont $w_{1},\dots,w_{N}$ et d’angles intérieurs $\lambda_{1}\pi,\dots,\lambda_{N}\pi$. Soit $\phi\colon \mathbb{H} \to P$ une uniformisante qui s’étend en un homéomorphisme sur chaque arête de $P$ et qui envoie l’infini sur le point~$w_{N}$. Il existe $N - 1$ nombres réels $\alpha_{1},\dots, \alpha_{N-1}$ tels que
\begin{equation}\label{eq:SC}
\phi(z) = C \int_{z_{0}}^{z}  \frac{dw}{(w-\alpha_{1})^{1-\lambda_{1}} \cdots (w-\alpha_{N-1})^{1-\lambda_{N-1}}}\,.
\end{equation}
Dans la preuve je vais utiliser cette formule avec des polygones spéciaux que j'introduis maintenant. Un {\em peigne} est une demi-bande horizontale avec des coupures horizontale. Plus précisément, soient $M\in \NN^{\ast}$ et $g \leq M -1$ des entiers naturels, on définit le peigne $\Pi_{M}(q_{1},\dots,q_{g};h_{1},\dots,h_{g})$ avec $q_{1}<\cdots<q_{g}$ des entiers contenus dans l'intervalle $[1,N-1]$ et $h_{i}>0$ par l'ensemble
\begin{equation}\label{eq:peigne}
\{ \xi\in \CC : \Re (\xi) >0 \text{ et } 0<\Im(\xi)<M\pi  \} \setminus \bigcup_{j=1}^{g} \{ \Im(\xi) =\pi q_{j} \text{ et } 0 < \Re(\xi) <h_{j} \}\,.
\end{equation}
La figure~\ref{fig:peigne} donne un exemple de peigne avec deux coupures horizontales, qui est lié à la différentielle de la figure~\ref{fig:diffcanong2}.

Ces rappels fait, je suis maintenant en mesure de donner la preuve du lemme~\ref{lm:existancecano}.
\begin{proof}
On considère sur une courbe totalement réelle $C_{D}$ une différentielle $\omega$  satisfaisant aux conditions de la proposition~\ref{prop:propdiffcancara}. Le quotient de $\omega$ par l'involution hyperelliptique donne une différentielle quadratique $\eta$ sur $\PP^{1}$ qui possède $g$ zéros doubles, $2g+2$ pôles simples et un unique pôle double dont le résidu quadratique est égal à~$1$. Une telle différentielle quadratique est représentée  dans la figure~\ref{fig:quotient} en prenant le quotient de la différentielle de la figure~\ref{fig:diffcanong2}. 
\begin{figure}[ht]
 \centering
\begin{tikzpicture}[scale=1,decoration={
    markings,
    mark=at position 0.5 with {\arrow[very thick]{>}}}]
   \fill[fill=black!10] (0,-3) coordinate (a0) -- (6,-3) coordinate (r2)-- (6,3) coordinate (r3) -- (0,3) coordinate (a0b) -- cycle;
   \draw (r2) -- (a0) -- (a0b) -- (r3);
    \draw[dotted] (r2) -- ++(.5,0);
\draw[dotted] (r3) -- ++(.5,0);

\draw[red,thick] (0,-2.02) coordinate(b0) -- (1,-2.02);
\draw[blue,thick] (0,-1.98) coordinate (a1) -- (1,-1.98) ;
\draw[green,thick]  (0,-1.02) coordinate (a2b) -- (1.2,-1.02);
\draw[black,thick] (0,-.98)  coordinate (b1b) -- (1.2,-.98);
\draw[red,thick] (0,2.02) coordinate (b0b) -- (1,2.02);
\draw[blue,thick] (0,1.98) coordinate(a1b) -- (1,1.98) ;
\draw[green,thick] (0,1.02)  coordinate (a2) -- (1.2,1.02);
\draw[black,thick] (0,.98) coordinate (b1) -- (1.2,.98);

\coordinate (b2) at (0,0);

\coordinate (c1b) at (1,-2); 
\coordinate (c2b) at (1.2,-1);
\coordinate (c3b) at (1.2,1); 
\coordinate (c4b) at (1,2);

\filldraw[fill=blue] (c1b) circle (2pt);
\filldraw[fill=blue] (c4b) circle (2pt);
\filldraw[fill=green] (c2b) circle (2pt);
\filldraw[fill=green] (c3b) circle (2pt);

\filldraw[fill=yellow] (a0)++(.1,0) arc [start angle=0, end angle=90, radius=1mm] -- ++(0,-.1);
\filldraw[fill=yellow] (a0b)++(.1,0)  arc [start angle=0, end angle=-90, radius=1mm]  -- ++(0,.1);
\filldraw[fill=red] (b0)++(.1,0)  arc [start angle=0, end angle=-90, radius=1mm]-- ++(0,.1);
\filldraw[fill=red] (b0b)++(.1,0)  arc  [start angle=0, end angle=90, radius=1mm] -- ++(0,-.1);
\filldraw[fill=purple] (a1)++(.1,0) arc [start angle=0, end angle=90, radius=1mm] -- ++(0,-.1);
\filldraw[fill=purple] (a1b)++(.1,0)  arc [start angle=0, end angle=-90, radius=1mm]  -- ++(0,.1);
\filldraw[fill=teal] (b1)++(.1,0)  arc [start angle=0, end angle=-90, radius=1mm]-- ++(0,.1);
\filldraw[fill=teal]  (b1b)++(.1,0)  arc  [start angle=0, end angle=90, radius=1mm] -- ++(0,-.1);
\fill (a2)++(.1,0) arc [start angle=0, end angle=90, radius=1mm]-- ++(0,-.1);
\fill (a2b)++(.1,0)  arc [start angle=0, end angle=-90, radius=1mm] -- ++(0,.1);
\filldraw[fill=white]  (b2) ++(0,.1)  arc [start angle=90, end angle=-90, radius=1mm];

\node[rotate=90] at (0,-2.5) {$1$};
\node[rotate=-90] at (0,2.5) {$1$};
\node[rotate=90] at (0,-1.5) {$2$};
\node[rotate=-90] at (0,1.5) {$2$};
\node[rotate=90] at (0,-.5) {$3$};
\node[rotate=-90] at (0,.5) {$3$};

\draw [decorate,decoration={brace}]	(0,2.07) -- (1,2.07) node [midway, above] {$\lambda_{1}$};
\draw [decorate,decoration={brace}]	(0,1.07) -- (1.2,1.07) node [midway, above] {$\lambda_{2}$};
\end{tikzpicture}
 \caption{La différentielle quadratique (multipliée par $\tfrac{r}{\pi}$) sur $\PP^{1}$ obtenue par quotient de la différentielle de la figure~\ref{fig:diffcanong2}. Les arêtes indexées par le même numéro sont identifiées deux à deux par une rotation d'angle~$\pi$ et une translation verticale.} \label{fig:quotient}
\end{figure}
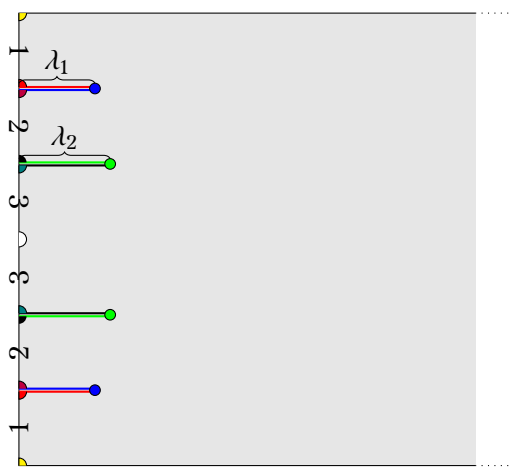
\par
 On peut alors restreindre $\eta$ au demi-plan supérieur~$\mathbb{H}$ et considérer son représentation plate obtenue en intégrant l'une des deux racines de~$\eta$. Cela revient, quitte à considérer l'autre racine, à restreindre la différentielle quadratique de la figure~\ref{fig:quotient} à la moitié inférieure du dessin. On obtient alors le peigne représenté dans la figure~\ref{fig:peigne} avec $h_{i} = \lambda_{i}/\pi$. Il est clair qu'en général on obtient un peigne, que je noterai~$\Pi_{\omega}$. Plus précisément, le peigne~$\Pi_{\omega}$ est donné par $\Pi_{r}(q_{1},\dots,q_{g};r\lambda_{1}/\pi,\dots,r\lambda_{g}/\pi)$, où $ q_{i} = \sum_{j=1}^{i}r_{j-1}$ avec les~$r_{j}$ définis dans le paragraphe précédent la proposition~\ref{prop:propdiffcancara} et les $\lambda_{i}$ les demi-périodes plates de~$\omega$. On note ~$\phi$ le biholomorphisme entre le peigne $\Pi_{\omega}$ et $\mathbb{H}$ normalisé en imposant
 \[\phi(a_{0}) = 0 ,\, \phi(b_{g})= i\pi r \text{ et } \phi(\infty) = \infty \,.\]
 Le biholomorphisme $\phi$ est donné par l'inverse de l'application de Schwarz-Christoffel:
 \begin{equation}\label{eq:SCconcret}
 \phi\colon \mathbb{H} \to \Pi_{\omega}: z \mapsto   r \int_{a_{0}}^{z}  \frac{(w-c_{1}) \cdots (w-c_{g})}{(w-a_{0})^{1/2}(w-b_{0})^{1/2}  \cdots (w-a_{g})^{1/2} (w-b_{g})^{1/2}}dw\,,
 \end{equation}
avec $a_{i},b_{j},c_{k}\in \RR$. On en déduit que  le polynôme $D$ est égal à
\[ D(w) = \prod_{i=0}^{g} (w-a_{i})(w-b_{i})\]
et que le polynôme canonique $R$ est
\[ R(w) = (w-c_{1}) \cdots (w-c_{g}) \,. \]
Maintenant, on définit sur la droite réelle la fonction
\begin{equation}
 P(x) = \cos (in \phi(x)) \,.
\end{equation}
Notons que cette fonction est analytique et s'étend à $\CC$ en une fonction holomorphe par la même formule. De plus, pour tout $i\in\left\{ 0,\dots,g \right\}$ elle possède $r_{i}$ zéros sur l'intervalle~$E_{i}$ et ne s'annule pas ailleurs. Comme $\phi$ est un homéomorphisme sur ces intervalles, les zéros de~$P$ sont simples.  On en déduit que $P$ est un polynôme de degré~$r$.
\par
Il reste à vérifier que $P$ vérifie l'équation de Pell-Abel sur $C_{D}$. Il est clair que la fonction $S(x)=\sin (in \phi(x))$ satisfaisait à l'équation $P^{2}+ S^{2} =1$. On en déduit que $S$ est un polynôme. Il reste à vérifier que $S$ est divisible par $D$. Cela est clairement impliqué par le fait que l'ensemble des racines de $S$ contient les points $a_{i}$ et~$b_{j}$.
\end{proof}

Avant de passer à la preuve  du théorème~\ref{thm:principal} je souhaite faire quelques commentaires sur le résultat que l'on vient de prouver.

L'un des points fondamentaux de la preuve est d'écrire $P$ comme le cosinus d'une intégrale hyperelliptique. Il existe d'autres moyens d'obtenir ce résultat. L'un des plus élégant est donné via les équations différentielles ordinaires. On pourra se reporter à la section~2 de \cite{bogarepr} pour obtenir cette version.

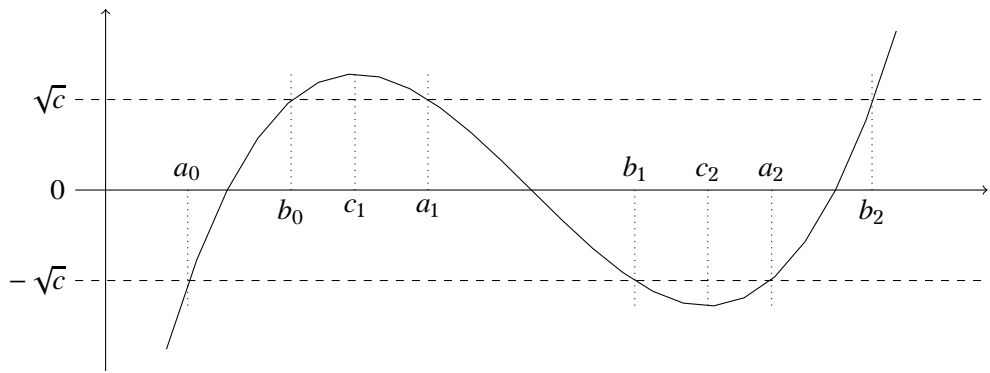
\begin{figure}[ht]
 \centering
\begin{tikzpicture}[scale=4,decoration={
    markings,
    mark=at position 0.5 with {\arrow[very thick]{>}}}]
  
  \draw[->] (.5,0) node[left] {$0$} -- (3.5,0);
  \draw[->] (.6,-.6) -- (.6,.6);
  \draw[dashed] (.5,.3) node[left] {$\sqrt{c}$} -- (3.5,.3);
   \draw[dashed] (.5,-.3) node[left] {$-\sqrt{c}$} -- (3.5,-.3);
\draw[domain=0.8:3.2] plot (\x,{(\x-1)*(\x-2)*(\x-3)});
\coordinate (a0) at (0.87,0);\node[above] at (a0) {$a_{0}$};
\draw[dotted] (a0) -- ++ (0,-.4);
\coordinate (b1) at (2.34,0);\node[above] at (b1) {$b_{1}$};
\draw[dotted] (b1) -- ++ (0,-.4);
\coordinate (a2) at (2.79,0);\node[above] at (a2) {$a_{2}$};
\draw[dotted] (a2) -- ++ (0,-.4);
\coordinate (b0) at (1.21,0);\node[below] at (b0) {$b_{0}$};
\draw[dotted] (b0) -- ++ (0,.4);
\coordinate (a1) at (1.66,0);\node[below] at (a1) {$a_{1}$};
\draw[dotted] (a1) -- ++ (0,.4);
\coordinate (b2) at (3.12,0);\node[below] at (b2) {$b_{2}$};
\draw[dotted] (b2) -- ++ (0,.4);

\coordinate (c1) at (1.42,0);\node[below] at (c1) {$c_{1}$};
\draw[dotted] (c1) -- ++ (0,.4);
\coordinate (c2) at (2.58,0);\node[above] at (c2) {$c_{2}$};
\draw[dotted] (c2) -- ++ (0,-.4);
\end{tikzpicture}
 \caption{Le polynôme $P$ associé à la différentielle canonique de la figure~\ref{fig:diffcanong2}.} \label{fig:polyasso}
\end{figure}

La figure~\ref{fig:polyasso} donne le polynôme associé à la différentielle de la figure~\ref{fig:diffcanong2}.
Le polynôme~$R$ est égal à  $(x-c_{1})(x-c_{2})$. De plus, $c_{i}$ est l'image de $z_{i}$ et $\bar{z}_{i}$ par l'involution hyperelliptique. La période plate $\lambda_{i}$ est simplement égale à l'intégrale de la forme $R/\sqrt{D} dx$ du point $b_{i-1}$ à $c_{i}$. Le fait que chaque intervalle~$E_{i}$ possède un unique zéro de $P$ est équivalent au fait que $r_{i}=1$. 
\par
On peut se demander si, étant donnée une différentielle canonique, on peut trouver de manière effective l'équation de la courbe hyperelliptique $C_{D}$ et le polynôme~$P$. Ce problème est très difficile en général, mais des résultats intéressants ont étés obtenus par \cite{akhi} et développés dans \cite{bogaeff} via les fonctions automorphes.
\par 
Il me semble enfin utile d'écrire explicitement l'un des points essentiel de la preuve du lemme~\ref{lm:existancecano}. Ce résultat formalise la relation entre le peigne de la figure~\ref{fig:peigne} et la différentielle de la figure~\ref{fig:diffcanong2}.
\begin{prop}\label{prop:peignediff}
 Soit $\Pi$ un peigne, alors la structure plate associée à la différentielle canonique associée à $\Pi$ est obtenue de la façon suivante.
 \begin{enumerate}
  \item Faire une réflexion de $\Pi$ le long de $\RR$. 
  \item Faire une rotation centrale au point $(0,0)$.
  \item Identifier par translation verticale les droites du bord du domaine ainsi créé.
  \item  Pour chaque segment obtenu par les images des dents du peigne, identifier par translation la partie supérieure du segment avec la partie inférieure du segment d'ordonné opposée.
 \end{enumerate}
\end{prop}
\smallskip
\par
Je donne maintenant la preuve du théorème~\ref{thm:principal}. Rappelons que celui-ci donne l'existence d'une solution primitive de degré $r$ de l'équation de Pell-Abel sur une courbe hyperelliptique totalement réelle de genre~$g$ si et seulement si $r\geq g+1$.
\begin{proof}
 Au vu de la remarque~\ref{rem:nonex}, il suffit de montrer qu'il existe des solutions de degré $r$ pour tout $r\geq g+1$. Donc par le lemme~\ref{lm:existancecano}, il suffit de montrer qu'il existe une différentielle canonique de degré $r$ sur une surface de genre~$g$ qui n'est pas de degré $r'$ pour un $r'<r$. Pour cela il suffit de considérer les entiers $r_{i} = 1$ pour $i\in \{0,\dots,g-1\}$ et $r_{g}=r-g$. On considère alors pour $q_{i}= \sum_{j=1}^{i}r_{j-1}$ le peigne 
 \[ \Pi_{r}(q_{1},\dots,g_{g};r\lambda_{1}/\pi,\dots,r\lambda_{g}/\pi)\]
 avec $\lambda_{i}>0$ pour tout $i\in\{1,\dots,g\}$. On obtient alors une différentielle satisfaisant les conditions de la proposition~\ref{prop:propdiffcancara} grâce à la proposition~\ref{prop:peignediff}.  Par le lemme~\ref{lm:existancecano}, cette différentielle est canonique d'ordre $r$ sur une surface totalement réelle de genre~$g$. Le choix des $r_{i}$ impliquent qu'elle est primitive. Donc la solution de l'équation de Pell-Abel associée à cette différentielle est primitive.
\end{proof}

On peut de plus déduire une description de l'ensemble des composantes connexes des lieux formés par les couples donnés par les différentielles canonique primitive de degré~$r$ sur les courbes hyperelliptiques totalement réelles de genre~$g$ sur lesquelles l'équation de Pell-Abel est résoluble. Un résultat similaire est donné par l'équation~(19) de \cite{bogaeff}.
\begin{cor}
 L'espace des modules des différentielles canonique primitive de degré~$r$ sur les courbes hyperelliptiques totalement réelles de genre~$g$ satisfaisant à l'équation de Pell-Abel est une union de cônes ouverts $(\mathbb{R}_{+}^{\ast})^{g-1}$, où cette union est prise sur toute les partitions  $(r_{0},\dots,r_{g})$ de $r$ à $g+1$ éléments tels que $\pgcd(r_{0},\dots,r_{g}) = 1$.
\end{cor}
\smallskip
\par
Pour terminer cette section, je donne une remarque sur le corps de définition des courbes considérées.
\begin{rem}
 On peut se demander si il existe des courbes algébriques définies sur $\mathbb{Q}$ ou un autre sous-corps de $\mathbb{R}$ avec des points de $r$-torsions (voir entre beaucoup d'autres références \cite{lepre}). Afin d'étudier ce problème, \cite{pakocombi} associe un graphe aux courbes hyperelliptiques avec un point de $r$-torsion. Ce graphe est la préimage du segment $[-\sqrt{c},\sqrt{c}]$ par le polynôme $P$. Dans le cas des courbes hyperelliptiques totalement réelles, on obtient une union de $g+1$ segments disjoints. Chacun de ces segments contient $r_{i}-1$ sommets en son intérieur. Donc ce cas ne semble pas être le cas le plus adéquat pour étudier ce problème.  
\end{rem}

\section{Les $k$-différentielles}
\label{sec:kdiff}

Dans cette section, je donne tout d'abord quelques rappels sur les $k$-différentielles méromorphes et la structure plate qui leur est associée pour tout $k\geq1$. Beaucoup plus d'informations peut se trouver dans \cite{mata,zorich,strebel} dans le cas des différentielles abéliennes et quadratiques et \cite{BCGGM3} dans celui des $k$-différentielles pour $k\geq3$. Une fois ces rappels fait, je donnerai la preuve du théorème~\ref{thm:proj}.

 \subsection{Rappels sur les $k$-différentielles}\label{sec:rappelsdiff}
Une {\em différentielle abélienne} est une paire $(C,\omega)$ où $C$ est une surface de Riemann et~$\omega$ est une section méromorphe non nulle du fibré cotangent $K$ de~$C$. Si la section est méromorphe et  ne possède que des pôles simples alors on dit qu'elle est de  {\em troisième espèce}. Dans le cas où $\omega$ est holomorphe, on dit quelle est de {\em première espèce}. 
\par
La différentielle~$\omega$ induit par intégration une structure de translation et une métrique plate sur la surface~$C$ privée des zéros  et des pôles de~$\omega$. Cette métrique s'étend aux zéros de $\omega$ en une métrique singulière. Plus précisément, un zéro d'ordre $m$ de $\omega$ correspond à une singularité conique d'angle $(m+1)2\pi$ de la métrique plate. De plus, un pôle simple de résidu $r$ correspond à un cylindre infini de circonférence $2i\pi r$. On peut donner une description au voisinage des pôles d'ordres supérieur ou égaux à $2$, toutefois on n'aura besoin que du cas des pôles doubles sans résidus. La structure plate au voisinage de tels pôles est isométrique au complémentaire d'un disque dans le plan muni de la métrique euclidienne standard. Pour le cas général, le lecteur pourra se reporter à \cite{strebel} pour plus de détails.   Réciproquement, on peut associer à une surface de translation satisfaisant certaines conditions une différentielle méromorphe. Grâce à cela, on peut effectivement vérifier que la surface de translation représentée dans la figure~\ref{fig:diffcanong2} correspond à une différentielle abélienne possédant les propriétés énoncées.
\par
Le théorème de Riemann-Roch implique qu'une différentielle $\omega$ possède $2g-2$ zéros comptés avec multiplicité. \'Etant donnée une partition $\mu = (m_{1},\dots,m_{n})$ de~$2g-2$, on définit donc l'espace des modules des différentielles abéliennes dont les ordres des zéros sont égaux à  $m_{1},\dots,m_{n}$. Cet espace des modules se nomme {\em strate} de différentielles abéliennes de type $\mu$ et se note $\omoduli(\mu)$. Ces strates sont des variétés orbifoldes complexes de dimension $2g-1+n$. 
\smallskip
\par
Il existe une compactification des strates de différentielles dont les éléments ont étés caractérisés dans \cite{BCGGM1}. Ces éléments sont  caractérisés par des objets, les {\em différentielles entrelacées}, qui sont donnés de la façon suivante.  C'est une collection de différentielles satisfaisant certaines conditions sur une surface de Riemann marquée stable. Ces différentielles sont obtenues en redimensionnant sur chaque composante de la limite les familles de différentielles. 
\smallskip
\par
La discussion précédente peut être étendue au cas des $k$-différentielles pour tout $k\geq2$. Une $k$-différentielle~$\xi$ est une section non nulle du produit tensoriel $k$ième du fibré canonique d'une surface de Riemann. Une $k$-différentielle en dite primitive si elle n'est pas la puissance d'une $k'$-différentielle avec $k'<k$. 

On obtient une surface plate en intégrant la racine $k$ième d'une $k$-différentielle~$\xi$.  Cette surface est une surface dont les changements de cartes, dans le complémentaires des singularités de $\xi$, sont donnés par une composition de translations et de rotations d'angles multiples de $2\pi/k$. Par exemple dans le cas d'une différentielle quadratique, ces changements de cartes sont des compositions de translations et de rotations d'angles~$\pi$. De manière équivalente, on peut obtenir une $k$-différentielle en identifiant les arêtes d'un polygone par des translations et des rotations d'angles multiples de $2\pi/k$. Une différentielle quadratique est donnée de la sorte sur la figure~\ref{fig:quotient}.

De plus, étant donnée une partition $\mu = (m_{1},\dots,m_{n})$ de~$k(2g-2)$, on considère les espaces des modules des $k$-différentielles dont les ordres des zéros et éventuellement des pôles sont égaux à  $m_{1},\dots,m_{n}$. Cet espace des modules se nomme {\em strate} de $k$-différentielles de type~$\mu$ et se note $\komoduli(\mu)$. Ces strates sont des variétés orbifoldes telles que les composantes paramétrant les $k$-différentielles primitives sont toutes de dimension complexe égale à $2g-2+n$.

 \subsection{Preuve du théorème~\ref{thm:proj}}

Je donne maintenant la preuve du théorème~\ref{thm:proj}. Rappelons que celui-ci donne l'existence d'une $k$-différentielle primitive avec un unique zéro sur une surface de Riemann hyperelliptique totalement réelle sauf dans le cas particulier des différentielles quadratiques sur les courbes de genre~$2$. 
\smallskip
\par
  Par le théorème~\ref{thm:principal}, pour tout $r\geq g+1$, il existe une courbe hyperelliptique totalement réelle~$C_{D}$ et des points $(z,z')$ conjugués par l'involution hyperelliptique tels que 
  \begin{equation}\label{eq:zerosn}
   rz - rz' \sim  \calO \,,
  \end{equation}
 où $\calO$ est le fibré trivial de $C_{D}$. De plus,  cette équation n'est satisfaite sur $C_{D}$ pour aucun $ m \in \left\{ 1,\dots,r-1 \right\}$.
De plus, comme ces points sont conjugués par l'involution hyperelliptique, on a l'égalité classique
  \begin{equation}\label{eq:cano}
   (g-1)z + (g-1)z' \sim K \,,
  \end{equation}
  où $K$ est le fibré canonique de $C_{D}$.
  Prenons maintenant $r=k(g-1)$, on obtient des équations~\eqref{eq:zerosn} et~\eqref{eq:cano}  que
  \begin{equation}
   2rz = kK \,.
  \end{equation}
  Donc $z$ est l'unique zéro d'une $k$-différentielle~$\xi$. Le fait que $r$ est minimal implique que~$\xi$ est une $k$-différentielle primitive. Comme $r\geq g+1$, on a l'existence d'une $k$-différentielle primitive qui possède un unique zéro pour tout $k\geq \tfrac{g+1}{g-1}$. L'unique restriction induite par cette inégalité est dans le cas $g=k=2$. Il est bien connu qu'il n'existe pas de différentielles quadratiques primitives en genre $2$ avec un unique zéro.   Cela conclut la preuve du théorème~\ref{thm:proj}.
 \smallskip
 \par
 Avant de passer à la section suivante, je souhaite noter les points suivants.
  \begin{rem}
  \begin{enumerate}
    \item Le lieu hyperelliptique est de dimension $2g-1$ et la projection de $\komoduli(k(2g-2))$ de dimension $2g-2$ dans l'espace des modules $\moduli$. Il n'est donc pas a priori évident que ces deux lieux s'intersectent.
   \item En genre $2$ on peut donner une représentation plate des $k$-différentielles de $\komoduli[2](2k)$ (voir par exemple \cite{CCM}). Je ne connais aucun moyen de savoir si cette représentation plate correspond à une $k$-différentielle sur une surface totalement réelle. De plus, dans le cas des courbes de genre $g\geq3$ je ne connais aucun moyen de déterminer si une représentation plate d'une telle $k$-différentielles est sur une courbe hyperelliptique (réelle ou même complexe). 
   \item Il serait intéressant de pouvoir donner une relation directe entre les représentations plates des différentielles canoniques et des $k$-différentielles induites.
 \end{enumerate}
 \end{rem}

\subsection{$k$-différentielles dont le diviseur est de support est conjugué}\label{sec:kdiffconjug}

Dans cette section, je vais généraliser le théorème~\ref{thm:proj} au cas des $k$-différentielles dont le diviseur est supporté sur des points $(z,\iota(z))$, où $\iota$ est l'involution hyperelliptique. Le résultat suivant donne une description de cette situation.
\begin{prop}\label{prop:kdiffconjug}
 Soient $C$ une courbe hyperelliptique de genre $g$ et $k\geq2$. S'il existe une $k$-différentielle  primitive~$\xi$ sur $C$ telle que $(\xi) = n z + n' \iota(z)$ avec $n>n'\geq 0$, alors l'équation de Pell-Abel est résoluble sur~$C$ et $2k(g-1)\geq n\geq (k+1)g - (k-1)$. 
 
 Réciproquement, si  $2k(g-1)\geq n\geq (k+1)g - (k-1)$ alors il existe une courbe hyperelliptique totalement réelle $C$ et une $k$-différentielle  primitive $\xi$ sur $C$ dont le diviseur est de la forme $(\xi) = n z + n' \iota(z)$.
\end{prop}
Notons que l'on retrouve le théorème~\ref{thm:proj} en posant $n=2k(g-1)$. De plus, nous ne considérons pas le cas $n=n'$ car il correspond à la puissance $k$ième d'une différentielle abélienne dont le diviseur est $(g-1)z + (g-1)\iota(z)$. La fin de cette section est constituée de la preuve de la proposition~\ref{prop:kdiffconjug}.
\smallskip
 \par
Soit $C$ une courbe hyperelliptique de genre $g$ telle qu'il existe une $k$-différentielle~$\xi$ holomorphe $\xi$ dont le diviseur est  $(\xi) = n z + n' \iota(z)$ avec $n>n'$. On supposera que $C$ est donnée par l'équation $y^{2} =\prod_{i=0}^{2g+1} (x-x_{i})$ et que les points $z$ et $\iota(z)$ sont les points à l'infini~$\infty_{\pm}$ de~$C$. Une base des $k$-différentielles holomorphes sur  $C$ est donnée par 
\[ x^{i} \frac{dx^{k}}{y^{k}}, \text{ pour }i\in \{0,\dots kg-k\} \text{ et } x^{i} \frac{dx^{k}}{y^{k-1}}, \text{ pour }i\in \{0,\dots kg-k-g-1\} \,. \]
 Donc une $k$-différentielle $\xi$ sur $C$ est de la forme
 \[\xi = \left(\sum_{i=0}^{kg-k} (\alpha_{j}x^{j} + \beta_{j}x^{j}y) \right) \frac{dx^{k}}{y^{k}} \,\]
où $\alpha_{j},\beta_{j}\in \CC $ et $\beta_{j}= 0$ pour $j\geq kg-k-g$.  De plus, le diviseur de la fonction $y$ sur~$C$ est 
\[(y) = \sum_{i=0}^{2g+1} (x_{i},0) - (g+1)(\infty_{-}+\infty_{+}) \]
et celui de la forme $dx$ est 
\[(dx) = \sum_{i=0}^{2g+1} (x_{i},0) - 2(\infty_{-}+\infty_{+}) \,.\]
On en déduit que  le support du diviseur de $dx^{k}/y^{k}$ est contenu dans les points $\infty_{\pm}$. Donc pour que le support du diviseur de $\xi$ soit contenu dans $\infty_{\pm}$, il faut que la fonction $\sum (\alpha_{j}x^{j} + \beta_{j}x^{j}y)$ soit un élément inversible de l'algèbre affine de~$C$. Comme la discussion de la section~\ref{sec:pelldebut}, s'étend sans difficultés aux courbes hyperelliptiques générales, on déduit de l'équation~\eqref{eq:Pell} que cet élément est une solution de l'équation de Pell-Abel.  Supposons qu'il existe une solution primitive de degré $r$ de Pell-Abel, alors les équations~\eqref{eq:zerosn} et~\eqref{eq:cano} donnent 
\begin{equation}
 (k(g-1) +r) \infty_{+} + (k(g-1)-r)\infty_{-} \sim kK \,.
\end{equation}
On a alors  $n=k(g-1) +r$ et le théorème~\ref{thm:principal} donne $r \geq g+1$, ce qui implique que $n \geq (k+1)g - (k-1)$. L'inégalité $2k(g-1)\geq n$ est une conséquence directe du fait que nous ne considérons que des $k$-différentielles holomorphes.
\par
La réciproque est une conséquence directe du théorème~\ref{thm:principal} et des calculs que l'on vient de réaliser.

 \section{Les points de torsion}
 \label{sec:tor}
 
Le but de cette section est de démontrer le théorème~\ref{thm:torsion} qui donne l'existence de points de $n$-torsion pour $n\geq 2g+1$.  Je donne quelques rappels avant de procéder à la preuve de ce résultat.

 \subsection{Rappels sur les points de torsion}
Soit $C$ une courbe hyperelliptique de genre $g\geq2$.  On considère un point de \Weierstrass~$W$ sur $C$ et l'application d'Abel-Jacobi $\varphi_{W}\colon C \to J_{C}$ qui envoie $W$  sur l'origine. Un {\em point de $n$-torsion} modulo~$W$ sur $C$  est un point $P$  tel que $\varphi_{W}(P)$ est un point de $n$-torsion de la jacobienne, mais pas de $n'$-torsion pour tout $0< n' < n$. Notons que les points de $2$-torsion modulo $W$ sont exactement les  autres points de \Weierstrass de~$C$.

On peut étendre ces considérations au cas des surfaces de Riemann de genre~$1$. Dans ce cas, on fixe un point arbitraire $W$ sur la surface $C$ de genre~$1$. Un {\em point de $N$-torsion} modulo~$W$ sur $C$  est un point $P$  tel que $P$ est un point de $N$-torsion sur la courbe elliptique $(C,W)$ mais pas de $n'$-torsion pour tout $0< n' < n$. J’appelle {\em points de \Weierstrass} les points de $2$-torsion de~$(C,W)$.

Pour terminer ces rappels, je souhaite insister sur le fait que l'ordre $n$ d'un point~$P$ dépend du choix du point de \Weierstrass. En effet,  prenons une surface de Riemann $C$ de genre~$1$ donnée par le quotient de $\CC$ par le réseau $\ZZ v_{1}\oplus \ZZ v_{2}$. Considérons les points~$W_{0}$ et $W_{1}$ de $C$ donnés respectivement par l'image de l'origine et de $\tfrac{1}{2} v_{1}$. Le point $P_{1}$ donné par l'image de $\tfrac{1}{6}v_{1}$ est de $6$-torsion modulo $W_{0}$ et de $3$-torsion modulo $W_{1}$. Cet exemple peut se généraliser au cas des courbes hyperelliptiques de genre $g\geq2$. 

 \subsection{Preuve du théorème~\ref{thm:torsion}}\label{sec:thmtorsion}

Je commence par montrer que pour $3\leq r \leq 2g$ il n'existe pas de points de $r$-torsion sur les courbes hyperelliptiques. Ce résultat est donné dans \cite[théorème~2.8]{Zarhin} et \cite[lemme 3.139]{drraabel}.  Je donne ici une preuve un peu différente, mais essentiellement équivalente à celles proposées par ces auteurs.
\par 
Supposons qu'il existe un tel point $P$ de $r$-torsion sur une courbe hyperelliptique~$C$. Par définition il existe un point de \Weierstrass~$W$ tel que $r(P-W)\sim \mathcal{O}$, où $\calO$ est le fibré trivial de~$C$. De plus, on a l'égalité classique $(2g-2)W \sim K$ avec $K$ le fibré canonique de~$C$. Donc l’existence de ce point de $r$-torsion est équivalente à l’existence d'une solution à l'équation 
$$rP + (2g-2-r)W \sim K\,.$$ 
Si $r\leq 2g-2$, la non existence de telle solution est une conséquence directe du fait que les diviseurs des différentielles sur les courbes hyperelliptiques sont invariants par l'involution hyperelliptique~$\iota$. Dans le cas $r=2g-1$, cela est une conséquence du théorème des résidus. En effet, il n'existe pas de différentielles avec un unique pôle simple sur une surface de Riemann compacte. Enfin, dans le cas $r=2g$, on utilise l'égalité classique $-2W + P +\iota(P) \sim \calO$ pour obtenir l'équation équivalente $(2g+1)P -\iota(P) \sim K$ qui n'a pas de solutions pour la même raison que dans le cas précédent.
\smallskip
\par
On montre maintenant que si $r\geq 2g+1$, alors il existe une courbe hyperelliptique totalement réelle de genre $g$ avec un point de $r$-torsion modulo un point de \Weierstrass. 
 Soit $W$ un point de \Weierstrass d'une courbe hyperelliptique totalement réelle $C_{D}$ qui possède une solution primitive d'ordre $r$ à l'équation de Pell-Abel. En notant $z$ et $z'$  les pôles $\infty_+$ et $\infty_-$ de la différentielle canonique associée, on déduit de $z + z' \sim 2W$ l'équation
 \begin{equation}\label{eq:torsion}
   2rz \sim 2r W \,.
 \end{equation}
Cela implique que le point $z$ est soit un point  $\ell$-torsion modulo $W$ avec $\ell | 2r$. Le fait que la solution est primitive de degré $r$ et que $z$ n'est pas de \Weierstrass implique que  soit $\ell= 2r$ ou $\ell = r$. Le reste de la preuve consiste à montrer que, quitte à choisir un autre point de \Weierstrass, le point $z$ est un point de $r$-torsion. On commence par traiter le cas du genre~$1$, puis le cas général par récurrence.
\smallskip
\par
Fixons un point $W_{0}$ sur une surface de Riemann $C$ de genre~$1$. Dans ce cas, il suffit de considérer un point primitif de $r$-torsion $z$ sur la courbe elliptique $(C,W_{0})$ pour obtenir le résultat. 
\smallskip
\par
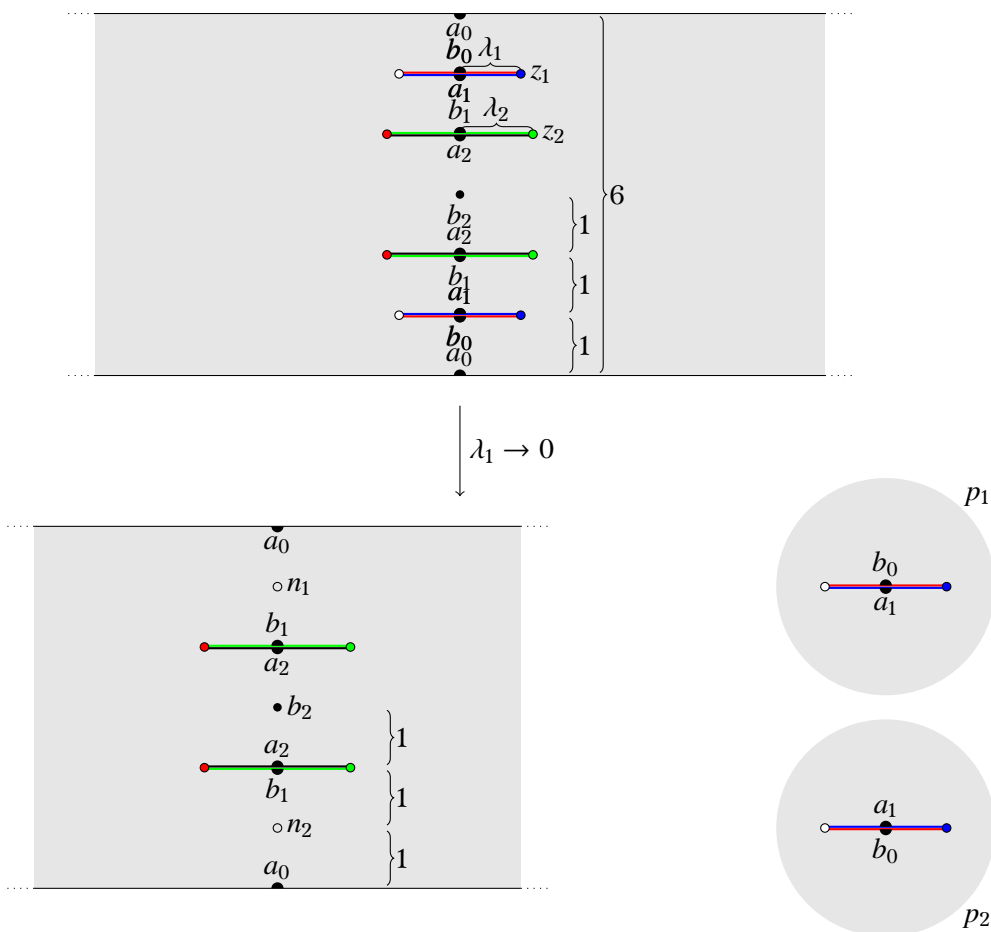
\begin{figure}[ht]
 \centering
\begin{tikzpicture}[scale=.8,decoration={
    markings,
    mark=at position 0.5 with {\arrow[very thick]{>}}}]
   \fill[fill=black!10] (-6,-3) coordinate (r1) -- (6,-3) coordinate (r2) coordinate[pos=.5] (a0) -- (6,3) coordinate (r3) -- (-6,3) coordinate (r4)coordinate[pos=.5] (a0b) -- cycle;
   \draw (r1) -- (r2);
   \draw (r3) -- (r4);
   \draw[dotted] (r1) -- ++(-.5,0);
    \draw[dotted] (r2) -- ++(.5,0);
\draw[dotted] (r3) -- ++(.5,0);
\draw[dotted] (r4) -- ++(-.5,0);

\draw[red,thick] (-1,-2.02) -- (1,-2.02) coordinate[pos=.5] (b0);
\draw[blue,thick] (-1,-1.98) -- (1,-1.98) coordinate[pos=.5] (a1);
\draw[green,thick] (-1.2,-1.02) -- (1.2,-1.02) coordinate[pos=.5] (b1);
\draw[black,thick] (-1.2,-.98) -- (1.2,-.98) coordinate[pos=.5] (a2);
\draw[red,thick] (-1,2.02) -- (1,2.02) coordinate[pos=.5] (b0b);
\draw[blue,thick] (-1,1.98) -- (1,1.98) coordinate[pos=.5] (a1b);
\draw[green,thick] (-1.2,1.02) -- (1.2,1.02) coordinate[pos=.5] (b1b);
\draw[black,thick] (-1.2,.98) -- (1.2,.98) coordinate[pos=.5] (a2b);

\coordinate (b2) at (0,0);

\coordinate (c1) at (-1,-2); 
\coordinate (c2) at (-1.2,-1);
\coordinate (c3) at (-1.2,1); 
\coordinate (c4) at (-1,2);
\coordinate (c1b) at (1,-2); 
\coordinate (c2b) at (1.2,-1);
\coordinate (c3b) at (1.2,1); 
\coordinate (c4b) at (1,2);

\filldraw[fill=white] (c1) circle (2pt);
\filldraw[fill=white] (c4) circle (2pt);
\filldraw[fill=red] (c2) circle (2pt);
\filldraw[fill=red] (c3) circle (2pt);

\filldraw[fill=blue] (c1b) circle (2pt);
\filldraw[fill=blue] (c4b) circle (2pt);
\filldraw[fill=green] (c2b) circle (2pt);
\filldraw[fill=green] (c3b) circle (2pt);

\node[right] at (c3b) {$z_{2}$};
\node[right] at (c4b) {$z_{1}$};

\foreach \i in {0,1}{
\fill (a\i)++(.1,0) arc [start angle=0, end angle=180, radius=1mm];
\node[above] at (a\i) {$a_{\i}$}; 
\fill (a\i b)++(-.1,0)  arc [start angle=-180, end angle=0, radius=1mm];
\node[below] at (a\i b) {$a_{\i}$}; 
\fill (b\i)++(-.1,0)  arc [start angle=-180, end angle=0, radius=1mm];
\node[below] at (b\i) {$b_{\i}$}; 
\fill (b\i b)++(.1,0)  arc  [start angle=0, end angle=180, radius=1mm];
\node[above] at (b\i b) {$b_{\i}$}; 
}
\fill (a2)++(.1,0) arc [start angle=0, end angle=180, radius=1mm];
\node[above] at (a2) {$a_{2}$}; 
\fill (a2b)++(-.1,0)  arc [start angle=-180, end angle=0, radius=1mm];
\node[below] at (a2b) {$a_{2}$}; 
\fill (b2) circle (2pt);\node[below] at (b2) {$b_{2}$};

\draw [decorate,decoration={brace}]	(1.8,-2.05) -- (1.8,-2.95) node [midway, right] {$1$};
\draw [decorate,decoration={brace}]	(1.8,-1.05) -- (1.8,-1.95) node [midway, right] {$1$};	
\draw [decorate,decoration={brace}]	(1.8,-.05) -- (1.8,-.95) node [midway, right] {$1$};	

\draw [decorate,decoration={brace}]	(2.3,2.95) -- (2.3,-2.95) node [midway, right] {$6$};	

\draw [decorate,decoration={brace}]	(0,2.07) -- (1,2.07) node [midway, above] {$\lambda_{1}$};
\draw [decorate,decoration={brace}]	(0,1.07) -- (1.2,1.07) node [midway, above] {$\lambda_{2}$};

\draw[->] (0,-3.5) -- (0,-5)  node [midway, right] {$\lambda_{1} \to 0$};	

\begin{scope}[yshift=-8.5cm,xshift=-3cm]
    \fill[fill=black!10] (-4,-3) coordinate (r1) -- (4,-3) coordinate (r2) coordinate[pos=.5] (a0) -- (4,3) coordinate (r3) -- (-4,3) coordinate (r4)coordinate[pos=.5] (a0b) -- cycle;
   \draw (r1) -- (r2);
   \draw (r3) -- (r4);
   \draw[dotted] (r1) -- ++(-.5,0);
    \draw[dotted] (r2) -- ++(.5,0);
\draw[dotted] (r3) -- ++(.5,0);
\draw[dotted] (r4) -- ++(-.5,0);

\draw[green,thick] (-1.2,-1.02) -- (1.2,-1.02) coordinate[pos=.5] (b1);
\draw[black,thick] (-1.2,-.98) -- (1.2,-.98) coordinate[pos=.5] (a2);
\draw[green,thick] (-1.2,1.02) -- (1.2,1.02) coordinate[pos=.5] (b1b);
\draw[black,thick] (-1.2,.98) -- (1.2,.98) coordinate[pos=.5] (a2b);

\coordinate (b2) at (0,0);

\coordinate (c1) at (-1,-2); 
\coordinate (c2) at (-1.2,-1);
\coordinate (c3) at (-1.2,1); 
\coordinate (c4) at (-1,2);
\coordinate (c1b) at (1,-2); 
\coordinate (c2b) at (1.2,-1);
\coordinate (c3b) at (1.2,1); 
\coordinate (c4b) at (1,2);

\filldraw[fill=red] (c2) circle (2pt);
\filldraw[fill=red] (c3) circle (2pt);

\filldraw[fill=green] (c2b) circle (2pt);
\filldraw[fill=green] (c3b) circle (2pt);

\foreach \i in {0,1}{
\fill (a\i)++(.1,0) arc [start angle=0, end angle=180, radius=1mm];
\node[above] at (a\i) {$a_{\i}$}; 
\fill (a\i b)++(-.1,0)  arc [start angle=-180, end angle=0, radius=1mm];
\node[below] at (a\i b) {$a_{\i}$}; 
\fill (b\i)++(-.1,0)  arc [start angle=-180, end angle=0, radius=1mm];
\node[below] at (b\i) {$b_{\i}$}; 
\fill (b\i b)++(.1,0)  arc  [start angle=0, end angle=180, radius=1mm];
\node[above] at (b\i b) {$b_{\i}$}; 
}
\fill (a2)++(.1,0) arc [start angle=0, end angle=180, radius=1mm];
\node[above] at (a2) {$a_{2}$}; 
\fill (a2b)++(-.1,0)  arc [start angle=-180, end angle=0, radius=1mm];
\node[below] at (a2b) {$a_{2}$}; 
\fill (b2) circle (2pt);\node[right] at (b2) {$b_{2}$};

\draw [decorate,decoration={brace}]	(1.8,-2.05) -- (1.8,-2.95) node [midway, right] {$1$};
\draw [decorate,decoration={brace}]	(1.8,-1.05) -- (1.8,-1.95) node [midway, right] {$1$};	
\draw [decorate,decoration={brace}]	(1.8,-.05) -- (1.8,-.95) node [midway, right] {$1$};	

\draw (0,2) circle (2pt);\node[right] at (0,2) {$n_{1}$}; 
\draw (0,-2) circle (2pt);\node[right] at (0,-2) {$n_{2}$}; 


      \fill[fill=black!10] (10,2) circle (1.8cm);
      \draw[red,thick] (9,2.02) -- (11,2.02) coordinate[pos=.5] (b0b);
\draw[blue,thick] (9,1.98) -- (11,1.98) coordinate[pos=.5] (a1b);

      \fill[fill=black!10] (10,-2) circle (1.8cm);
            \draw[red,thick] (9,-2.02) -- (11,-2.02) coordinate[pos=.5] (b0);
\draw[blue,thick] (9,-1.98) -- (11,-1.98) coordinate[pos=.5] (a1);

\coordinate (c1) at (9,-2); 
\coordinate (c2) at (9,2);
\coordinate (c1b) at (11,2); 
\coordinate (c2b) at (11,-2);

\filldraw[fill=white] (c1) circle (2pt);
\filldraw[fill=white] (c2) circle (2pt);

\filldraw[fill=blue] (c1b) circle (2pt);
\filldraw[fill=blue] (c2b) circle (2pt);

\fill (a1)++(.1,0) arc [start angle=0, end angle=180, radius=1mm];
\node[above] at (a1) {$a_{1}$}; 
\fill (a1b)++(-.1,0)  arc [start angle=-180, end angle=0, radius=1mm];
\node[below] at (a1b) {$a_{1}$}; 
\fill (b0)++(-.1,0)  arc [start angle=-180, end angle=0, radius=1mm];
\node[below] at (b0) {$b_{0}$}; 
\fill (b0b)++(.1,0)  arc  [start angle=0, end angle=180, radius=1mm];
\node[above] at (b0b) {$b_{0}$}; 

\node at (11.5,3.5) {$p_{1}$};
\node at (11.5,-3.5) {$p_{2}$};
\end{scope}

\end{tikzpicture}
 \caption{La différentielle canonique (multipliée par $\tfrac{r}{\pi}$) de la figure~\ref{fig:diffcanong2} et la différentielle entrelacée obtenue par la  dégénération $\lambda_{1} \to 0$. La différentielle en bas à droite est obtenue en divisant cette famille de différentielles par $\lambda_{1}$.} \label{fig:diffcanong2degen}
\end{figure}

Je traite maintenant le cas des courbes de genre $g\geq2$. On se donne une différentielle canonique $\eta$ de degré $r$ sur une courbe hyperelliptique totalement réelle de genre~$g$. On peut alors déformer cette différentielle de la façon suivante, représentée sur la figure~\ref{fig:diffcanong2degen}. 
Rappelons qu'une période plate  de $\omega$ est l'intégrale de $\omega$ entre deux zéros conjugués par l'involution hyperelliptique. On fait tendre une période plate, disons $\lambda_{1}$, vers $0$ en maintenant les autres périodes plates constantes. Notons que par le lemme~\ref{lm:existancecano}, on obtient une famille de différentielles canoniques primitives de degré~$r$. 

La différentielle entrelacée (voir la section~\ref{sec:rappelsdiff}) obtenue à la limite  est définie sur deux composantes qui forment une courbe banane. Une composante est une courbe hyperelliptique de genre $g-1$ et l'autre composante est une courbe de genre~$0$. Rappelons qu'une courbe banane est une courbe composée de deux composantes irréductibles lisses s'intersectant en deux n\oe uds. Dans la figure~\ref{fig:diffcanong2degen}, cette courbe banane est obtenue en collant les points $n_{i}$ aux pôles $p_{i}$. La différentielle $\eta_{1}$ sur la composante~$C_{1}$ de genre $g-1$ est la différentielle que l'on voit à la limite. Elle appartient à la strate $\omoduli[g-1](1,\dots,1,-1,-1)$ avec $2g-2$ zéros simples. Pour obtenir la différentielle~$\eta_{2}$, on multiplie la famille de différentielles par~$1/\lambda_{1}$ et faisons tendre $\lambda_{1}$ vers~$0$. On obtient alors la différentielle représentée à droite de la figure~\ref{fig:diffcanong2degen}. Cette différentielle est définie sur la sphère de Riemann, appartient à la strate $\omoduli[0](1,1,-2,-2)$ et les résidus aux pôles sont nuls.  Notons que par le lemme~\ref{lm:existancecano} la différentielle $\eta_{1}$ est une différentielle canonique de degré $r$ sur la courbe hyperelliptique totalement réelle~$C_{1}$ de genre $g-1$. De plus, on peut choisir les $r_{i}$ de telle façon que la différentielle $\eta_{1}$ est primitive. Par récurrence,  pour tout $r\geq 2g+1$ il existe  un point de \Weierstrass~$W_{1}$ sur~$C_{1}$ tel que la différence les pôles sont de $r$-torsion modulo le point $W_{1}$. La preuve du théorème~\ref{thm:torsion} se conclut en utilisant le fait que  l'ordre d'un point de $r$-torsion qui reste fini par déformation est constant.

 \addcontentsline{toc}{section}{Références}
\bibliographystyle{alpha}
\bibliography{../biblio}

\bigskip
\noindent
\small{Quentin Gendron\\
Centro de Investigaci\'on en Matem\'aticas,\\ Guanjuato, Gto.,\\
AP 402, CP 36000,\\
México\\
{\em courriel:} quentin.gendron@cimat.mx} 
 \end{document}